\documentclass[10pt,a4paper]{article}
\usepackage{amsfonts}
\usepackage{amsmath, amssymb}
\usepackage{makeidx}
\makeindex

\newtheorem{theorem}{Theorem}[section]
\newtheorem{lemma}[theorem]{Lemma}

\newtheorem{remark}[theorem]{Remark}

\newcommand{\R}{{{\mathbb R}}}

\numberwithin{equation}{section}

\begin{document}
\author{Yanjin Wang \thanks{Corresponding Email: wangyj@ms.u-tokyo.ac.jp}\thanks{The author wishes to express his deep
gratitude to Prof. Hitoshi Kitada for his constant encouragement,
kind guidance and helpful discussions. The author also thanks Dr.
Xiaoqun Zhang, Laboratoire LMAM of Universit\'e de Bretagne-Sud, for
her help. The study is supported by Japanese Government
Scholarship.}\\{\small \textit{Graduate School of Mathematics,
University of Tokyo, 3-8-1 Komaba,}}
\\{\small\textit{ Meguro, Tokyo, 153-8914, Japan}}}

\title{Finite time blow-up results for the damped wave equations with arbitrary initial energy in an inhomogeneous medium}
\date{}
\maketitle

\begin{center}
\textbf{Abstract}
\end{center}
In this paper we consider the long time behavior of solutions of
the initial value problem for the damped wave equation of the form
\begin{eqnarray*}
 u_{tt}-\rho(x)^{-1}\Delta u+u_t+m^2u=f(u)
\end{eqnarray*}
with some $\rho(x)$ and $f(u)$ on the whole space $\R^n$ ($n\geq
3$).

For the low initial energy case, which is the non-positive initial
energy, based on concavity argument we prove the blow up result. As
for the high initial energy case, we give out sufficient conditions
of the initial datum such that the corresponding solution blows up
in finite time.

\begin{center}
\textbf{R\'esum\'e}
\end{center}
Dans cet article, nous consid\'erons le comportement asymptotique
des solutions d'un probl\'me \'a valeur initiale pour l\'equation
d'ondes att\'enu\'ees  suivante
 \begin{eqnarray*}
  u_{tt}-\rho(x)^{-1}\Delta u+u_t+m^2u=f(u)
 \end{eqnarray*}
o\'u $\rho(x)$ et $f(u)$ sur l'espace $\R^n$ ($n\geq3$).

Dans le cas d'\'energie initiale faible, qui est l'\'energie
initiale non positive, en nous basant  sur un argument de
concavit\'e,  nous establissons un  r\'esultat d'explosion. Quant au
cas d'\'energie initiale forte, nous donnons  des conditions
suffisantes sur la donn\'ee initiale pour que l'explosion de la
solution  correspondante ait lieu en un temps fini.



\section{Introduction}
In this paper our aim is to study a class of wave equations in the
following form
\begin{eqnarray}
\left\{\begin{array}{r@{,\quad  }l}
 u_{tt}-\rho(x)^{-1}\Delta u+u_t+m^2u=f(u)& (t,x)\in [0,T)\times \R^n\\
 u(0,x)=u_0(x)& x\in\R^n\\
 u_t(0,x)=u_1(x)&x\in\R^n
\end{array}
\right.\label{DWE}
\end{eqnarray}
where $\Delta$ is Laplacian operator on $\R^n$ ($n\geq 3$),
$u_0(x)$ and $u_1(x)$ are real valued functions, $m$ is a real
constant (the case, $m=0$, is called as the mass free case; $m\ne
0$ as the mass
case), $\rho(x)$ satisfies the following condition\\
\textbf{(H)} \textit{$\rho(x)>0$ for every $x\in\R^n$, $\rho\in
\mathcal{C}^{0,\gamma}(\R^n)$ with $\gamma\in(0,1)$, and $\rho\in
L^{n/2}(\R^n)\cap L^\infty(\R^n)$.}

The wave equations (\ref{DWE}) appear in applications in various
areas of mathematical physics
\cite{01Antman}\cite{02reed}\cite{03Zauderer}, as well as in
geophysics and ocean acoustics, where, for example, the
coefficient $\rho(x)$ represents the speed of sound at the point
$x\in\R^N$\cite{04Klibanov} , in other words, $\rho(x)\ne
constant$ implies that the medium is inhomogeneous, where the
sound travels.

For the nonlinear power, throughout the paper we make the following
assumption: the nonlinear power $f(s)$ satisfies that there exists
some constant $\epsilon>0$ such that
\begin{eqnarray}
f(s)s\geq (2+\epsilon)F(s)\label{PowerAssum}
\end{eqnarray}
for any $s\in\R$, where
\begin{eqnarray}
F(\zeta)=\int_0^\zeta f(\kappa)d\kappa. \label{F-Def}
\end{eqnarray}
For such nonlinear power, it was first stated for abstract wave
equations with $\rho(x)=1$ by Levine \cite{05Levine}. And then
Cazenave \cite{07Cazenave} also considered it for Klein-Gordon
equations.

Before going any further, we first briefly introduce some research
works for the wave equation (\ref{DWE}) with
$\rho(x)=\mathrm{constant}\ne 0$ (without loss of generality let
$\rho(x)=1$), obviously it does not satisfy the assumption (H).
For the general nonlinear power $f(u)$ with (\ref{PowerAssum}) it
was firstly considered for some abstract wave equations in
\cite{05Levine}, where Levine proved the blow up result when the
initial energy was negative. But mostly the results about Cauchy
problem for the wave equation were investigated the typical form
of nonlinear power as
\begin{eqnarray}
f(u)=|u|^{p-1}u \label{Power-Special}
\end{eqnarray}
where $1<p<\frac{n+2}{n-2}$. Here we note that the above power
satisfies the condition (\ref{PowerAssum}).  For the power
(\ref{Power-Special}), the wave equations with damping term were
studied by many authors.
 It is well known that the local solution blows up in finite time
 when the initial energy is negative.
 For global existence and nonexistence of solutions for
Cauchy problem of the equation (\ref{DWE}) with $\rho(x)=1$,
(\ref{Power-Special}) and (possibly nonlinear) damping term, we here
refer to
\cite{08Georgiev}\cite{09Ikehata}\cite{10Ikehata}\cite{06Levine}\cite{11Levine}\cite{16Nakao}\cite{12Pucci}\cite{13Vitillaro}.
In special, recently the wave equation with damping term was
considered in \cite{14Levine}, where Levine and Todorova showed that
for arbitrarily positive initial energy there are choices of initial
datum such that the local solution blows up in finite time.
Subsequent Todorova and Vitillaro \cite{23Todorova} established more
precise result about the existence of initial values such that the
corresponding solution blows up in finite time for arbitrarily high
initial energy. More recently, Gazzola and Squassina
\cite{15Gazzola} established sufficient conditions of initial datum
with arbitrarily positive initial energy such that the corresponding
solution  blows up in finite time for the wave equation with linear
damping and (\ref{Power-Special}) in the mass free case on a bounded
Lipschitz subset of $\R^n$. And the author \cite{28Wang} establish
the blow up result with arbitrarily positive initial energy for some
Klein-Gordon equations on the whole space $\R^n$.

For the case $\rho(x)=1$ and $m=0$, we also note that it was
considered that  $f(u)=|u|^p$ by some authors. Here we just refer
to the papers \cite{17Todorova}\cite{18Ikehata} and some papers
cited therein.

Now we return to the equation (\ref{DWE}) with some general
$\rho(x)$. For the linear case, $f(u)=0$, Eidus \cite{24Eidus} first
studied the existence of solutions for linear wave equation
(\ref{DWE}). Then Karachalios and Stavrakakis \cite{19Karachalis}
studied the existence of the solution of the damped wave equation
(\ref{DWE}) with some nonlinear power. And they \cite{20Karachalis}
also established the results about the global existence and blow up
of solutions for the equation (\ref{DWE}) with (H) and
(\ref{Power-Special}) in the free mass case by potential wall
method, which was firstly developed by Sattinger \cite{22Sattinger}.
Their blow up result was under the condition that the initial energy
was negative. Recently, for the equation (\ref{DWE}) with (H) and
(\ref{Power-Special}) Zhou \cite{21Zhou} investigated the global
existence and blow up result including the
 mass free case and mass case. Zhou established the blow up result
when the initial energy was less than a positive constant, which
dependent on the function $\rho(x)$ (H). But in all the above works,
the high initial energy case was not considered  for the equation
(\ref{DWE}) under the assumption (H). Moreover, there is no one who
consider the equation (\ref{DWE}) with some general nonlinear power,
for example, like (\ref{PowerAssum}).

The main purpose of this paper is to establish the blow up result
for the equation (\ref{DWE}) with (H) and (\ref{PowerAssum}) when
the initial energy is high. Based on a concavity argument, which was
originally introduced by Levine \cite{05Levine}\cite{06Levine}, we
first establish the blow up result when the initial energy is
non-positive. Thus we extend the blow up results \cite{20Karachalis}
and \cite{21Zhou}. Note that, when the initial energy vanishes, the
blow up result was also established in \cite{21Zhou}, where it
needed the assumption $\rho(x)\in L^1(\R^n)$ with (H) and $\int
\rho(x)u_0(x)u_1(x)dx\geq 0$. Here in some sense (See Remark 2.4) we
improve the blow up result in \cite{21Zhou}. As for the arbitrarily
positive initial energy case, we establish the blow up result under
some conditions of $(u_0, u_1)$ on the whole space $\R^n$.  To the
best of our knowledge this is the first blow up result with high
initial energy for the equation (\ref{DWE}) with (H). The work is
motivated by \cite{15Gazzola}. But Our proof is different with
\cite{15Gazzola}, and very simple.
 Additionally, our proof is also valid for the case $\rho(x)=constant\ne 0$.
So at last we also make some remarks on the case
$\rho(x)=constant\ne 0$. We find that, if $\rho(x)$ satisfies (H)
then the mass $m$ does not affect the blow up result, but if
$\rho=constant\ne 0$ it will affect the blow up result, that is, if
$m=0$ and $\rho(x)=1$ then the blow up result is obtained only on a
bounded subset of $\R^n$. Indeed, by (H) we see that $\rho(x)$ will
rapidly enough decreasing at infinity, thus it make us possibly
consider the equation (\ref{DWE}) on the whole space $\R^n$ in the
mass free case.

The paper is composed of three sections. In the next section we
will denote some notations, and state our main results. The last
section is the main part, we prove the blow up results there. In
addition, some remakes are made on the case of $\rho(x)=1$.

\section{Principal Result}
In order to state our main results, we briefly mention here some
facts, notations and known results. We denote by $\|\cdot\|_q$ the
$L^q(\R^n)$ norm for $1\leq q\leq \infty$, and we define the
spaces: $H^1(\R^n)=\{u\in L^2(\R^n);\
\|u\|_{H^1(\R^n)}=\|(1-\Delta)^{1/2}u\|<\infty\}$, and
$H_0^1(\R^n)=\{u\in H^1(\R^n);\ \mathrm{supp}(u)$ $\mathrm{is\
compact\ in\ }\R^n\}$. For simplicity we will denote $\int_{\R^n}$
by $\int$. The notation $t\rightarrow T^-$ means $t< T$ and
$t\rightarrow T$.

As \cite{25Kozono}, we introduce the function space
$X^{1,2}(R^n)$, which is defined as the closure of $C_0^\infty$
functions with respect to the energy norm $\|u\|_{X^{1,2}}:=\int
|\nabla u(x)|^2dx$, that is,
\begin{eqnarray}
X^{1,2}(\R^n)=\{u\in L^{\frac{2n}{n-2}}(\R^n) : \nabla u\in
(L^2(\R^n))^n\}.
\end{eqnarray}
And it is known that $X^{1,2}$ is embedding continuously in
$L^{\frac{2n}{n-2}}$, which means that there exists a constant $k>0$
such that
\begin{eqnarray}
\|u\|_{\frac{2n}{n-2}}\leq k\|u\|_{X^{1,2}}.\label{EmbT}
\end{eqnarray}
Following (\ref{EmbT}), we have the following inequality
\cite{26Brown}.
\begin{lemma}
 Suppose $\rho\in L^{\frac{n}{2}}(\R^n)$ and $n\geq 3$. Then there
exists a constant $\alpha>0$ such that
\begin{eqnarray}
\int |\nabla u(x)|^2dx\geq \alpha\int \rho(x)|u(x)|^2dx
\end{eqnarray}
for every $u\in C_0^\infty(\R^n)$. Moreover,
$\alpha=k^{-2}\|\rho\|^{-1}_{\frac{n}{2}}$ where $k$ is defined in
(\ref{EmbT}).
\end{lemma}
In addition, the weighted space $L^2_{\rho}(\R^n)$ is defined to be
the closure of $C_0^\infty(\R^n)$ functions with respect to the
inner product
\begin{eqnarray}
(u, v)_{L^2_\rho}:=\int \rho(x)u(x)v(x)dx,
\end{eqnarray}
and norm
\begin{eqnarray}
\|u\|_{L^2_\rho}^2&=&(u, u)_{L^2_\rho}.
\end{eqnarray}
For local existence of solutions of the equation (\ref{DWE}), we
state
\begin{theorem}
Under the assumption (H). Let the initial datum $(u_0, u_1)\in
X^{1,2}(\R^n)\times L_\rho^2(\R^n)$, and $f$ satisfying the
following conditions: $f(0)=0$ and
\begin{eqnarray}
|f(\lambda_1)-f(\lambda_2)|\leq c
(|\lambda_1|^{p-1}+|\lambda_2|^{p-1})|\lambda_1-\lambda_2|
\end{eqnarray}
for all $\lambda_1, \lambda_2\in\R$, some constant $c>0$, and
\begin{eqnarray}
1<p<\frac{n}{n-2} \mathrm{\ when\ } n\geq 3.
\end{eqnarray}
Then there exists a unique local solution $u(t,x)$ of the equation
(\ref{DWE}) on a maximal time interval $[0, T_{\max})$ satisfying
\begin{eqnarray}
u\in C([0, T_{\max}); X^{1,2}(\R^n)) \mathrm{\ and\ } u_t\in C([0, T_{\max}); L_\rho^2),\\
u(0, x)=u_0(x)\mathrm{\ and\ } u_t(0, x)=u_1(x).
\end{eqnarray}
In addition, $u(t,x)$ satisfies
\begin{eqnarray}
E(0)-E(t)=\int_0^t\|u_\tau(\tau,\cdot)\|_{L_\rho^2}^2d\tau\geq
0\label{EnergyDec}
\end{eqnarray}
for every $t\in [0, T_{\max})$, where
\begin{eqnarray}
E(t)=\frac{1}{2}\int\left(\rho(x)|u_t(t,x)|^2+|\nabla
u(t,x)|^2+m^2\rho(x)|u(t,x)|^2-2\rho(x)F(u(t,x))\right)dx,\nonumber
\\
\ \label{EnergyDef}
\end{eqnarray}
where $F(u)$ is defined in (\ref{F-Def}).
\end{theorem}

This result can be proved by Banach fixed point theorem. The proof
follows from the weighted-norm Lebesgue space of the corresponding
theorem for the wave equations of Kirchhoff type \cite{27Ono}.

If $T_{\max}<\infty$, then the local solution is said to blow up
in finite time $T_{\max}$. Otherwise, $T_{\max}=\infty$, the
corresponding local solution is global.

Next we state our first blow up result for the equation
(\ref{DWE}) with (H) and (\ref{PowerAssum}) in the non-positive
initial energy case.
\begin{theorem} Under the assumptions (H) and (\ref{PowerAssum}). If the nonzero initial datum
 $(u_0,u_1)\in X^{1,2}(\R^n)\times
L_\rho^2(\R^n)$satisfies
\begin{eqnarray}
E(0)<0,
\end{eqnarray}
or
\begin{eqnarray}
\int \rho(x)u_0(x)u_1(x)dx\geq 0 \mathrm{\ if\ }
E(0)=0,\label{ThmLE01}
\end{eqnarray}
then the corresponding local solution of the equation (\ref{DWE})
blows up in finite time $T_{\max}<\infty$, that is,
\begin{eqnarray}
\lim_{t\rightarrow T_{\max}^-}\int \rho(x)|u(t,x)|^2dx=\infty.
\end{eqnarray}
\end{theorem}
\begin{remark}
In the case $E(0)=0$, Zhou \cite{21Zhou} also established the blow
up result for the equation (\ref{DWE}) with (H) and
(\ref{Power-Special})
 under some another assumptions as $\int \rho(x)u_0(x)u_1(x)dx\geq 0$ and $\rho\in
 L^1(\R^n)$. But in the above theorem, we remove $\rho\in
 L^1(\R^n)$.
 Thus, in this sense we improve the result \cite{21Zhou}.
\end{remark}

To state our main blow up result for the arbitrarily positive
initial energy case, we introduce a function as follows
\begin{eqnarray}
I(u)=\int (|\nabla u(x)|^2+m^2\rho(x)|u(x)|-\rho(x)f(u(x))u(x))dx.
\end{eqnarray}

Now we introduce our main blow up result for the equation
(\ref{DWE}) in the arbitrarily positive initial energy case, as far
as we know, which is the first blow up result for the equation
(\ref{DWE}) with (H) on the whole space $\R^n$.
\begin{theorem}
Under the assumptions (H) and (\ref{PowerAssum}). If the initial
datum $(u_0,u_1)\in X^{1,2}(\R^n)\times L_\rho^2(\R^n)$satisfies
\begin{eqnarray}
&&E(0)>0,\label{ThmHE01}\\
&&I(u_0)<0,\label{ThmHE02}\\
&&\int \rho(x)u_0(x)u_1(x)dx\geq 0,\label{ThmHE03}\\
&&\|u_0\|_{L_\rho^2}^2>\frac{2(2+\epsilon)}{m^2\epsilon}E(0)
\mathrm{\ when\ } m\ne
0,\label{ThmHE04}\\
&&\|u_0\|_{L_\rho^2}^2>\frac{2(2+\epsilon)}{\min\{1,
\alpha\}\epsilon}E(0) \mathrm{\ when\ } m=0,\label{ThmHE05}
\end{eqnarray}
where $\epsilon$ and $\alpha$ are stated in (\ref{PowerAssum}) and
Lemma 2.1, respectively.  Then the corresponding local solution of
the equation (\ref{DWE}) blows up in finite time $T_{\max}<\infty$,
that is,
\begin{eqnarray}
\lim_{t\rightarrow T_{\max}^-}\int \rho(x)|u(t,x)|^2dx=\infty.
\end{eqnarray}
\end{theorem}
\begin{remark}
We note that, for the case $E(0)<0$, by (\ref{EnergyDef}) and
(\ref{EnergyDec}) it is valid that $I(u(t,\cdot))<0$ for every
$t\in[0, T_{\max})$.
\end{remark}
Reading Theorem 2.3, 2.5 and Remark 2.6, naturally one considers how
about the local solution when the initial data satisfies $E(0)>0$
and $I(u_0)>0$. Indeed, for this case, being similar as the argument
with $m=0$ and $f(u)=|u|^{p-1}u$ \cite{20Karachalis},  by a
potential wall method we can also obtain the global existence of
solutions of the equation (\ref{DWE}) with (H) and
(\ref{PowerAssum}) when the positive initial energy is small enough.
Here we omit it. Furthermore, it is still open that whether there
exists a global solution for wave equations when the initial energy
is arbitrarily high.

\section{Proof of the main theorems}
In this section, we prove Theorem 2.3 and 2.5 based on concavity
argument. Next we first claim two lemmas. The following lemma is
basic.
\begin{lemma}
Let $T>0$ and $H(t)$ be a Lipschitzian function over $[0, T)$.
Assume that $H(0)\geq 0$ and
\begin{eqnarray}
\frac{d}{dt}H(t)+H(t)>0
\end{eqnarray}
for every $t\in[0, T)$. Then $H(t)>0$ for every $t\in (0, T)$.
\end{lemma}
Following the way \cite{15Gazzola}, by Lemma 3.1 we can obtain next
lemma. For the convenient of readers and completeness of the paper
we here still give out a proof.
\begin{lemma}
Assume that $(u_0, u_1)\in X^{1,2}(\R^n)\times L^2_\rho(\R^n)$
satisfies that
\begin{eqnarray}
\int \rho(x)u_0(x)u_1(x)dx\geq 0.\label{u0u1>0}
\end{eqnarray}
If the corresponding local solution $(u(t,x), u_t(t,x))\in C([0,
T_{\max}), X^{1,2}(\R^n)\times L^2_\rho(\R^n))$ is such that
\begin{eqnarray}
I(u(t,\cdot))<0\label{I<0}
\end{eqnarray}
for every $t\in[0, T_{\max})$, then $\|u(t,\cdot)\|_{L^2_{\rho}}^2$
is strictly increasing on $[0, T_{\max})$.
\end{lemma}
\textbf{Proof.} Since $u(t,x)$ is the local solution of the equation
(\ref{DWE}), then we easily have
\begin{eqnarray}
\frac{1}{2}\frac{d^2}{dt^2}\int \rho(x)|u(t,x)|^2dx&=&\int
\rho(x)|u_t(t,x)|^2dx+
\int\rho(x)u(t,x)u_t(t,x)\nonumber\\
&&-\int\left(m^2\rho(x)|u(t,x)|^2+|\nabla
u(t,x)|^2-\rho(x)u(t,x)f(u(t,x))\right)dx,
\end{eqnarray}
for every $t\in[0, T_{\max})$.

As a result, it follows
\begin{eqnarray*}
\frac{d^2}{dt^2}\int \rho(x)|u(t,x)|^2dx+\frac{d}{dt}\int
\rho(x)|u(t,x)|^2dx=2(\|u(t,x)\|_{L^2_\rho}^2-I(u(t,\cdot)))
\end{eqnarray*}
By (\ref{I<0}) and the above equation we have
\begin{eqnarray}
\frac{d^2}{dt^2}\int \rho(x)|u(t,x)|^2dx+\frac{d}{dt}\int
\rho(x)|u(t,x)|^2dx>0 \label{UStric-Increas}
\end{eqnarray}
for every $t\in[0, T_{\max})$.

Here we let $\displaystyle H(t)=\frac{d}{dt}\int
\rho(x)|u(t,x)|^2dx$, then as \cite{20Karachalis} we see that the
function $H(t)$ is Lipschitzian function over $[0, T_{\max})$. Thus,
from Lemma 3.1 and (\ref{UStric-Increas}) it follows that
$\|u(t,x)\|^2_{L_\rho^2}$ is strictly increasing on $[0, T_{\max})$.
\begin{flushright}
$\Box$
\end{flushright}

\textbf{Proof of Theorem 2.3} We first define the following
auxiliary function
\begin{eqnarray}
G(t)=\int
\rho(x)|u(t,x)|^2dx+\int_0^t\|u(\tau,\cdot)\|_{L_\rho^2}^2d\tau+(T_0-t)\|u_0\|_{L_\rho^2}^2+\zeta(T_1+t)^2,\label{G-def}
\end{eqnarray}
where the constants, $T_0>0, T_1>0, \zeta>0$, will be determined
later.

We then have
\begin{eqnarray}
G^\prime(t)&=&\frac{d}{dt}G(t)\nonumber\\
&=&2\int \rho(x)u(t,x)u_t(t,\cdot)dx+
\|u(t,\cdot)\|_{L_\rho^2}^2-\|u_0\|_{L_\rho^2}^2+2\zeta(T_1+t)\nonumber\\
&=&2\int \rho(x)u(t,x)u_t(t,\cdot)dx+ 2\int_0^t(u(\tau,\cdot),
u_\tau(\tau,\cdot))_{L_\rho^2}d\tau +2\zeta(T_1+t),\label{G1diff}
\end{eqnarray}
and
\begin{eqnarray}
\frac{1}{2}G^{\prime\prime}(t)&=&\int
\rho(x)|u_t(t,x)|^2dx+\int\rho(t,x)u(t,x)u_{tt}(t,x)dx\nonumber\\
&&+\int\rho(x)u(\tau,x)u_\tau(t,x)dx+\zeta\nonumber\\
&=&\int \rho(x)|u_t(t,x)|^2dx+\int
\rho(x)u(t,x)f(u(t,x))dx\nonumber\\
&&-\int\left(|\nabla
u(t,x)|^2+m^2\rho(x)|u(t,x)|^2\right)dx+\zeta\label{G2diff}
\end{eqnarray}

\textbf{Case I: $E(0)<0$.} By (\ref{PowerAssum}),
(\ref{EnergyDec}) and (\ref{EnergyDef}), we see that
\begin{eqnarray}
\int \rho(x)F(u(t,x))dx&=& \frac{1}{2}\int
\left(m^2\rho(x)|u(t,x)|^2+\rho(x)|u_t(t,x)|^2+ |\nabla
u(t,x)|^2\right)dx\nonumber\\
&&-E(0)+\int_0^t\|u_\tau(\tau,x)\|_{L_\rho^2}^2d\tau\nonumber\\
&\leq&
\frac{1}{2+\epsilon}\int\rho(x)u(t,x)f(u(t,x))dx\label{PTh231}
\end{eqnarray}
for all $t\in[0, T_{\max})$.

Thus, from (\ref{G2diff}) it follows that
\begin{eqnarray}
G^{\prime\prime}(t)&\geq& (4+\epsilon)\int
\rho(x)|u_t(t,x)|^2dx+2(2+\epsilon)\int_0^t
\|u(\tau,\cdot)\|_{L^{2}_{\rho}}^2d\tau\nonumber\\
&&+\epsilon\int\left(|\nabla
u(t,x)|^2+m^2\rho(x)|u(t,x)|\right)dx-2(2+\epsilon)E(0)+2\zeta.
\label{G2diff-2}
\end{eqnarray}
We now let the constant $\zeta$ satisfy
\begin{eqnarray}
\displaystyle 0<\zeta\leq-2E(0).\label{zeta-Def}
\end{eqnarray}
Then it follows from (\ref{zeta-Def})
\begin{eqnarray}
-2(2+\epsilon)E(0)+2\zeta\geq (4+\epsilon)\zeta.\label{zeta1}
\end{eqnarray}

Obviously, $G^{\prime\prime}(t)>0$ on $[0, T_{\max})$. Moreover, we
can take $T_1>0$ large sufficiently such that
\begin{eqnarray}
G^\prime(0)=2\int\rho(x)u_0(x)u_1(x)+2BT_1>0,\label{thm23-1}
\end{eqnarray}
and
\begin{eqnarray}
\frac{\epsilon}{2}\left(\int \rho(x)u_0(x)u_1(x)dx+\zeta
T_1\right)>\int \rho(x)|u_0(x)|^2dx.\label{thm23-2}
\end{eqnarray}

Thus, by (\ref{thm23-1}) we see that $G(t)>0$, $G^\prime(t)>0$ and
$G^{\prime\prime}(t)>0$ for every $t\in [0, T_{\max})$. That is,
$G(t)$ and $G^\prime(t)$ is strictly increasing on $[0,T_{\max})$.
Then we let
\begin{eqnarray}
A&=&\int\rho(x)|u(t,x)|^2dx+\int_0^t\|u(\tau,\cdot)\|_{L^2_\rho}^2d\tau+\zeta(T_1+t)^2\label{A-def}\\
B&=&\frac{1}{2}G^{\prime}(t) \label{B-def}\\
C&=&\int\rho(x)|u_t(t,x)|^2dx+\int_0^t\|u_\tau(\tau,\cdot)\|_{L^2_\rho}^2d\tau+\zeta.\label{C-def}
\end{eqnarray}

For every $t\in [0, T_0]$ we obviously have
\begin{eqnarray}
G(t)\geq A,
\end{eqnarray}
and by (\ref{G2diff-2}) and (\ref{zeta1})
\begin{eqnarray}
G^{\prime\prime}(t)\geq (4+\epsilon)C.
\end{eqnarray}

 We now let $T_0$ sufficiently large and satisfy
\begin{eqnarray}
T_0\geq\frac{4G(0)}{\epsilon G^\prime(0)}.\label{T0-def}
\end{eqnarray}
Noting the inequalities (\ref{G-def}), (\ref{G1diff}) and
(\ref{thm23-2}), we see that the definition of $T_0$ as
(\ref{T0-def}) is reasonable.

And suppose the solution $u(t,x)$ exists on $[0, T_0]$. Then it
follows that
\begin{eqnarray}
G^{\prime\prime}G(t)-\frac{4+\epsilon}{4}(G^\prime(t))^2\geq
(4+\epsilon)(AC-B^2)
\end{eqnarray}
for every $t\in[0,T_0]$

By a simple computation we see that
\begin{eqnarray}
As^2-2Bs+C&=&\int
\rho(x)(su(t,x)+u_t(t,x))^2dx\nonumber\\
&&+\int_0^t\|su(\tau,\cdot)+u_\tau(\tau,\cdot)\|^2_{L^2_\rho}d\tau+\zeta
(s(T_1+t)+1)^2\nonumber\\
&\geq& 0
\end{eqnarray}
for every $s\in \R$ and $t\in [0, T_0]$, which means that
$(2B)^2-4AC\leq 0$.

Thus we see that
\begin{eqnarray}
G^{\prime\prime}G(t)-\frac{4+\epsilon}{4}(G^\prime(t))^2\geq 0
\end{eqnarray}

Since $\displaystyle\frac{4+\epsilon}{4}>1$, we put $\displaystyle
\alpha=\frac{\epsilon}{4}$. Then we have
\begin{eqnarray}
\frac{d}{dt}G^{-\alpha}(t)&=&-\alpha G^{-\alpha-1}G^\prime(t)<0\\
\frac{d^2}{dt^2}G^{-\alpha}(t)&=&-\alpha
G^{-\alpha-2}\left[G^{\prime\prime}G(t)-\frac{4+\epsilon}{4}(G^\prime(t))^2\right]\nonumber\\
&\leq&0\label{Gconcave}
\end{eqnarray}
for every $t\in[0, T_0]$, which means that the function
$G^{-\alpha}$ is concave. Obviously $G(0)>0$, then from
(\ref{Gconcave}) it follows that the function
$G^{-\alpha}\rightarrow 0$ when $t<T_{\max}$ and $\displaystyle
t\rightarrow T_{\max}$ ($\displaystyle T_{\max}<\frac{4
G(0)}{\epsilon G^\prime(0)}\leq T_0$). Noting the assumption that
the solution exists on $[0,T_0]$, where $T_0$ is defined as
(\ref{T0-def}), thus we see that there exists a finite time
$T_{\max}>0$ such that
\begin{eqnarray}
\lim_{t\rightarrow T_{\max}^-}\|u(t,\cdot)\|_{L^2_\rho}^2=\infty,
\end{eqnarray}
which implies that the corresponding solution $u(t,x)$ of the
equation (\ref{DWE}) blows up in finite time $T_{\max}<\infty$.

\textbf{Case II: $E(0)=0$ and $\int \rho(x)u_0(x)u_1(x)dx\geq 0$.}
By (\ref{EnergyDec}) and (\ref{EnergyDef}) we have
\begin{eqnarray}
\int \left(|\nabla u(t,x)|^2+m^2\rho(x)|u(t,x)|^2\right)dx- 2\int
\rho(x)F(u(t,x))dx\leq 0
\end{eqnarray}
for every $t\in[0, T_{\max})$.

And noting the fact that $\displaystyle\int\rho(x)F(u(t,x))dx\ne 0$
we obtain by (\ref{PowerAssum})
\begin{eqnarray}
2\int \rho(x)F(u(t,x))dx< \int \rho(x)f(u(t,x))u(t,x)dx.
\end{eqnarray}
We then get
\begin{eqnarray}
I(u(t,x))<0
\end{eqnarray}
for every $t\in[0, T_{\max})$.

Thus by (\ref{ThmLE01}) and Lemma 3.2 we see that
$\|u(t,\cdot)\|_{L_\rho^2}^2$ is strictly increasing on $[0,
T_{\max})$.

In this case we still use the auxiliary function $G(t)$ as
(\ref{G-def}).

Thus, according to the proof of Case I, by (\ref{ThmLE01}) we see
that $G(t)>0$, $G^\prime(t)>0$, $G^{\prime\prime}(t)>0$ on $(0,
T_{\max})$, that is to say, $G(t)$ and $G^{\prime}(t)$ is strictly
increasing over $[0, T_{\max})$.

And as (\ref{G2diff-2}) we also have
\begin{eqnarray}
G^{\prime\prime}(t)&\geq& (4+\epsilon)\int
\rho(x)|u_t(t,x)|^2dx+2(2+\epsilon)\int_0^t
\|u(\tau,\cdot)\|_{L^{2}_{\rho}}^2d\tau\nonumber\\
&&+\epsilon\int\left(|\nabla
u(t,x)|^2+m^2\rho(x)|u(t,x)|^2\right)dx+2\zeta \nonumber\\
&\geq& (4+\epsilon)\int \rho(x)|u_t(t,x)|^2dx+2(2+\epsilon)\int_0^t
\|u(\tau,\cdot)\|_{L^{2}_{\rho}}^2d\tau\nonumber\\
&&+\epsilon(m^2+\alpha)\int\rho(x)|u_0(x)|^2dx+2\zeta,
\label{G2diff-3}
\end{eqnarray}
where the last inequality comes from Lemma 2.1 and Lemma 3.2.

Now we let the constant $\zeta$ satisfy
\begin{eqnarray}
0<\zeta\leq \frac{\epsilon}{2+\epsilon}\left(\alpha
+m^2\right)\|u_0\|^2_{L_\rho^2}
\end{eqnarray}
for the mass free case or the mass case, and the other positive
constants, $T_0$ and $T_1$, be large such that
\begin{eqnarray}
&&T_0\geq\frac{4G(0)}{\epsilon G^\prime(0)},\\
&&\frac{\epsilon}{2}\left(\int \rho(x)u_0(x)u_1(x)dx+\zeta
T_1\right)>\int \rho(x)|u_0(x)|^2dx.
\end{eqnarray}

Then by the same argument as Case I, we can claim that the
corresponding local solution of the equation (\ref{DWE}) blows up in
finite time.

Thus the proof of Theorem 2.3 is completed.
\begin{flushright}
$\Box$
\end{flushright}

In the following part we will process Theorem 2.5. Next lemma is the
crux to prove Theorem 2.5.
\begin{lemma}
Under the assumptions on $\rho(x)$, $f(u)$ and $(u_0,u_1)$ in
Theorem 2.6, then the corresponding local solution $(u(t,x),
u_t(t,x))\in C([0, T_{\max}), X^{1,2}(\R^n)\times L^2_\rho(\R^n))$
satisfies
\begin{eqnarray}
&&I(u(t,\cdot))<0, \label{I(t)<0}\\
&&\|u(t,\cdot)\|_{L_\rho^2}^2> \frac{2(2+\epsilon)}{m^2\epsilon}E(0)
\mathrm{\ when\ } m\ne 0,\label{Lemma33-1}\\
&&\|u(t,\cdot)\|_{L_\rho^2}^2>
\frac{2(2+\epsilon)}{\min\{1,\alpha\}\epsilon}E(0) \mathrm{\ when\ }
m= 0,\label{Lemma33-2}
\end{eqnarray}
for every $t\in[0, T_{\max})$.
\end{lemma}
\textbf{Proof.} Here the proof is by a contradiction argument. We
assume that (\ref{I(t)<0}) is not true over $[0, T_{max})$, that is,
there exists a time $T>0$ such that
\begin{eqnarray}
T=\min\{t\in(0, T_{\max}); I(u(t,\cdot))=0\}.
\end{eqnarray}

\textbf{Case I: $m\ne 0$.}
 Since $I(u(t,\cdot))<0$ on $[0, T)$, by Lemma 3.2 we see that
$\|u(t,\cdot)\|_{L^2_\rho}^2$ is strictly increasing over $[0, T)$,
which implies that
\begin{eqnarray}
\|u(t,\cdot)\|_{L^2_\rho}^2>\|u_0\|_{L^2_\rho}^2>
\frac{2(2+\epsilon)}{m^2\epsilon}E(0),
\end{eqnarray}
for $m\ne 0$ and every $t\in(0, T)$.

And by the continuity of $\|u(t,\cdot)\|^2_{L^2_\rho}$ at $t$, we
see that
\begin{eqnarray}
\|u(T,\cdot)\|_{L^2_\rho}^2>
\frac{2(2+\epsilon)}{m^2\epsilon}E(0).\label{Contradiction1}
\end{eqnarray}

On the other hand, by (\ref{EnergyDec}) and (\ref{EnergyDef})
 we see that
\begin{eqnarray}
m^2\|u(T,\cdot)\|_{L^2_\rho}^2+\|\nabla u(T,\cdot)\|^2-2\int
\rho(x)F(u(T,x))dx\leq 2E(T)\leq 2E(0).\label{Ctd21}
\end{eqnarray}
Moreover, noting the assumption $I(u(T,\cdot))=0$ and
(\ref{PowerAssum}), we then have
\begin{eqnarray}
\|u(T,\cdot)\|_{L^2_\rho}^2+\|\nabla u(T,\cdot)\|^2\geq
(2+\epsilon)\int \rho(x)F(u(T,x))dx.\label{Ctd22}
\end{eqnarray}
Combining (\ref{Ctd21}) and (\ref{Ctd22}) we then obtain
\begin{eqnarray}
m^2\|u(T,\cdot)\|_{L^2_\rho}^2+\|\nabla u(T,\cdot)\|^2\leq
\frac{2(2+\epsilon)}{\epsilon}E(0).\label{Contradiction2}
\end{eqnarray}
Obviously there is a contradiction between (\ref{Contradiction1})
and (\ref{Contradiction2}). Thus we have proved that
\begin{eqnarray}
I(u(t,\cdot))<0.\label{ILeq0}
\end{eqnarray}
for every $t\in[0, T_{\max})$.

By Lemma 3.2 we see that $\|u(t,\cdot)\|_{L^2_\rho}^2$ is strictly
increasing on $t$ if $I(u(t,\cdot))<0$ for every $t\in [0,
T_{\max})$. Thus, (\ref{ILeq0}) implies that
\begin{eqnarray}
\|u(t,\cdot)\|_{L^2_\rho}^2>\frac{2(2+\epsilon)}{\epsilon}E(0)
\end{eqnarray}
for every $t\in [0, T_{\max})$.

Hereunto the proof for Case I, $m\ne 0$, is accomplished.

\textbf{Case II: $m=0$.} As the argument for
(\ref{Contradiction1}), we can also obtain
\begin{eqnarray}
\|u(T,\cdot)\|_{L^2_\rho}^2>
\frac{2(2+\epsilon)}{\min\{1,\alpha\}\epsilon}E(0).\label{Contradiction3}
\end{eqnarray}
Since $m=0$, then the inequality (\ref{Contradiction2}) is rewritten
as
\begin{eqnarray}
\|\nabla u(T,\cdot)\|^2\leq \frac{2(2+\epsilon)}{\epsilon}E(0).
\end{eqnarray}
By Lemma 2.1, we see that
\begin{eqnarray}
\alpha \|u(t,x)\|_{L^2_\rho}^2\leq
\frac{2(2+\epsilon)}{\epsilon}E(0).\label{Contradiction4}
\end{eqnarray}
Thus by (\ref{Contradiction3}) and (\ref{Contradiction4}) we obtain
that the assumption, $I(u(T,\cdot))=0$, is  wrong. That is to say,
it is valid that $I(u(t,\cdot))<0$ for every $t\in[0, T_{\max})$.

Similarly, we also get
\begin{eqnarray}
\|u(t,\cdot)\|_{L^2_\rho}^2>\frac{2(2+\epsilon)}{\min\{1,\alpha\}\epsilon}E(0)
\end{eqnarray}
for every $t\in [0, T_{\max})$.

Thus all the proof of Lemma 2.3 has been completed.
\begin{flushright}
$\Box$
\end{flushright}
 \textbf{Proof of Theorem 2.5.} Here we still use the auxiliary
 function $G$, defined as (\ref{G-def}). We have
\begin{eqnarray}
G^{\prime\prime}(t)&=&2\int
\rho(x)|u_t(t,x)|dx+2\int\rho(x)f(u(t,x))u(tx,dx)\nonumber\\
&&-2\int(|\nabla u(t,x)|^2+m^2\rho(x)|u(t,x)|^2)dx+2\zeta\nonumber\\
&=&2\int
\rho(x)|u_t(t,x)|dx-2I(u(t,\cdot))+2\zeta\label{Thm24G2diff}
\end{eqnarray}
By Lemma 3.3, we see that
\begin{eqnarray}
G^{\prime\prime}(t)>0
\end{eqnarray}
for every $t\in[0, T_{\max})$.

And from (\ref{ThmHE03}) we see that $G^\prime(t)>0$ for every
$t\in(0, T_{\max})$. Thus, it comes that $G(t)$ and $G^\prime(t)$ is
strictly increasing on $[0, T_{\max})$.

Obviously the inequality (\ref{PTh231}) is also valid here for every
$t\in[0, T_{\max})$. Then by Lemma 2.1, Lemma 3.2,
(\ref{Lemma33-1}), (\ref{Lemma33-2}) and (\ref{Thm24G2diff}) we have
\begin{eqnarray}
G^{\prime\prime}(t)&\geq& (4+\epsilon)\int
\rho(x)|u_t(t,x)|^2dx+2(2+\epsilon)\int_0^t
\|u_\tau(\tau,\cdot)\|_{L^{2}_{\rho}}^2d\tau\nonumber\\
&&+\epsilon\int\left(|\nabla
u(t,x)|^2+m^2\rho(x)|u(t,x)|^2\right)dx-2(2+\epsilon)E(0)+2\zeta\nonumber\\
&\geq& (4+\epsilon)\int \rho(x)|u_t(t,x)|^2dx+2(2+\epsilon)\int_0^t
\|u_\tau(\tau,\cdot)\|_{L^{2}_{\rho}}^2d\tau\nonumber\\
&&+\epsilon\int
(m^2+\alpha)\rho(x)|u_0(x)|^2dx-2(2+\epsilon)E(0)+2\zeta
\end{eqnarray}
for every $t\in[0, T_{\max})$.

By (\ref{ThmHE04}) and (\ref{ThmHE05}) we see that
$$\frac{\epsilon(m^2+\alpha)}{2+\epsilon}\int \rho(x)|u_0(x)|^2dx-2E(0)>0,$$
for every $m\in\R$.

 We now let $\zeta$ satisfy
\begin{eqnarray}
0<\zeta\leq \frac{\epsilon(m^2+\alpha)}{2+\epsilon}\int
\rho(x)|u_0(x)|^2dx-2E(0)
\end{eqnarray}
for every $m\in\R$, the other positive constants, $T_0$, and $T_1$,
be large such that
\begin{eqnarray}
&&T_0\geq\frac{4G(0)}{\epsilon G^\prime(0)},\\
&&\frac{\epsilon}{2}\left(\int \rho(x)u_0(x)u_1(x)dx+\zeta
T_1\right)>\int \rho(x)|u_0(x)|^2dx.
\end{eqnarray}

We next let $A, B, C$ denote the same terms as
(\ref{A-def}),(\ref{B-def}) and (\ref{C-def}), respectively.

And assume that the solution $u(t,x)$ exists on $[0, T_0]$. Then we
have
$$G(t)\geq A$$
and
$$G^{\prime\prime}(t)\geq (4+\epsilon)C$$
for
every $t\in[0, T_0]$.

Thus by the same way as Theorem 2.3 it comes that
\begin{eqnarray}
G(t)G^{\prime\prime}(t)-\frac{4+\epsilon}{4}(G^\prime(t))^2\geq
(4+\epsilon)(AC-B^2)\geq 0
\end{eqnarray}

As the proof of Theorem 2.3, by a concavity argument we can also
obtain that, there exists  finite time $T_{\max}<\infty$ such that
\begin{eqnarray}
\lim_{t\rightarrow T_{\max}^-}\|u(t,\cdot)\|_{L^2_{\rho}}^2=\infty.
\end{eqnarray}
which implies the corresponding solution $u(t,x)$ of the equation
(\ref{DWE}) blows up in finite time $T_{\max}<\infty$.
\begin{flushright}
$\Box$
\end{flushright}
\begin{remark}
Reading the proof of Theorem 2.5, we can easily use a similar way to
the case $\rho(x)=1$ and $m\ne 0$ on the function space
$H_0^1(\R^n)\times L^2(\R^n)$.
\end{remark}
\begin{remark}
In the mass free case, we see that in the proof of Theorem 2.5 it is
necessary to use Lemma 2.1, which may be called as a general
Poincar\'e inequality. But it is well-known that Poincar\'e
inequatlity is valid on a bounded set. That is to say, In the mass
free case with $\rho(x)=1$ we cannot obtain a blow up result as
Theorem 2.5 on the whole space $\R^n$ when the initial energy is
high. For this case, Gazzola and Squassina \cite{15Gazzola} have
established the blow up result with arbitrarily positive initial
energy on a bounded Lipschitz subset of $\R^n$.
\end{remark}

\bibliographystyle{amsplain}

\begin{thebibliography}{10}
\bibitem{01Antman} S. S. Antman, The equations for large vibrations of strings, Amer. Math. Monthly 87(1980) 359-370.

\bibitem{26Brown} K.J. Brown and N. Stavrakakis, Global bifurcation
results for a semilinear elliptic equation on all of $\R^n$, Duke
Math. J., 85(1996) 77-94.

\bibitem{24Eidus} D. Eidus, The Cauchy problem for the wave equation in an
inhomogeneous medium, Comm. Partial Differential Equations 20
(1995)1589-1603.

\bibitem{07Cazenave}T. Cazenave, Uniform estimates for solutions of nonlinear Klein-Gordon equations, J. Funct.
Anal. 60 (1985) 36-55.

\bibitem{15Gazzola}F. Gazzola and M. Squassina, Global solutions and finite
time blowup for damped semilinear wave equations, Ann. I.H.
Poincar\`e-AN 23(2006)185-207.

\bibitem{08Georgiev}V. Georgiev, G. Todorova, Existence of a solution of the wave equation with nonlinear damping
and source term, J. Differential Equations 109 (1994) 295-308.

\bibitem{09Ikehata} R. Ikehata, Some remarks on the wave equations with nonlinear
damping and source terms, Nonlinear Anal. 27 (1996) 1165-1175.

\bibitem{10Ikehata} R. Ikehata, T. Suzuki, Stable and unstable sets for evolution
equations of parabolic and hyperbolic type, Hiroshima Math. J. 26
(1996) 475-491.

\bibitem{18Ikehata}R. Ikehata and K. Tanizawa, Global existence of solutions
for semilinear damped wave equations in $\R^N$ with noncompactly
supported initial data. Nonlinear Anal. 61 (2005) 1189-1208.

\bibitem{19Karachalis} N.I. Karachalios and N.M. Stavrakakis, Existence of a global attractor
for semilinear dissipative wave equations on $R\sp N$. J.
Differential Equations 157 (1999)183-205.

\bibitem{20Karachalis} N.I. Karachalios and N.M. Stavrakakis, Global existence and blow-up results
for some nonlinear wave equations on $\R^N$. Adv. Differential
Equations 6 (2001) 155-174.

\bibitem{04Klibanov}M. V. Klibanov, Global convexity in a three-dimensional inverse acoustic problem, SIAM J.
Math. Anal. 28(1997) 1371-1388.

\bibitem{25Kozono} H. Kozono and H. Sohr, New a prior estimates for
the Stokes equations in exterior domains, Indiana Univ. Math. J.
40(1991) 1-27.

\bibitem{05Levine}H. A. Levine, Instability and nonexistence of global solutions to nonlinear wave equations of
the form Putt = .Au+F(u), Trans. Amer. Math. Soc. 192(1974) 1-21.

\bibitem{06Levine}H.A. Levine, Some additional remarks on the nonexistence of global
solutions to nonlinear wave equations, SIAM J. Math. Anal. 5
(1974) 138-146.

\bibitem{11Levine} H.A. Levine, J. Serrin, Global nonexistence theorems for
quasilinear evolution equations with dissipation, Arch. Rational
Mech. Anal. 137 (1997) 341-361.

\bibitem{14Levine}H.A. Levine, G. Todorova, Blow up of solutions of the Cauchy
problem for a wave equation with nonlinear damping and source
terms and positive initial energy, Proc. Amer. Math. Soc. 129
(2001) 793-805.

\bibitem{16Nakao} M. Nakao and K. Ono,
Existence of global solutions to the Cauchy problem for the
semilinear dissipative wave equations, Math. Z. 214 (1993)
325-342.

\bibitem{27Ono} K. Ono, On global existence , asymptotic stability
and blowing up of solutions for some degenerate nonlinear wave
equations of Kirchhoff type with a strong dissipation, Mathematical
Methods in the Applied Science, 20(1997)151-177.

\bibitem{12Pucci} P. Pucci, J. Serrin, Global nonexistence for abstract
evolution equations with positive initial energy, J. Differential
Equations 150 (1998) 203-214.

\bibitem{02reed}M. Reed and B. Simon, Methods of Modern Mathematical Physics. III. Scattering Theory (Academic
Press, Harcourt Brace Jovanovich, Publishers, New York¨CLondon,
1979).

\bibitem{22Sattinger}D.H. Sattinger, On global solution of nonlinear
hyperbolic equations, Arch. Rational Mech. Anal. 30(1968)148-172.

\bibitem{17Todorova} G. Todorova and B. Yordanov,
Critical exponent for a nonlinear wave equation with damping. J.
Differential Equations 174 (2001) 464-489.

\bibitem{23Todorova}G. Todorova and E. Vitillaro,
Blow-up for nonlinear dissipative wave equations in$\R^n$, J.
Math. Anal. Appl. 303 (2005)242-257.

\bibitem{13Vitillaro} E. Vitillaro, Global existence theorems for a class of
evolution equations with dissipation, Arch. Rational Mech. Anal.
149 (1999) 155-182.

\bibitem{28Wang}Y. Wang, A sufficient condition for finite time
blow up of the Klein-Gordon equation with arbitrarily positive
initial energy. Submitting for publication.

\bibitem{03Zauderer} E. Zauderer, Partial Differential Equations of Applied Mathematics, Second edition, Pure and
Applied Mathematics, A Wiley-Intersci. Publ. (John Wiley $\&$
Sons, Inc., New York, 1989).

\bibitem{21Zhou} Y. Zhou, Global existence and nonexistence for a nonlinear wave equation
with damping and source terms. Math. Nachr. 278 (2005) 1341-1358.


\end{thebibliography}

\end{document}